\documentclass{amsart}
\usepackage[utf8]{inputenc}

\usepackage[english]{babel}

\usepackage{amssymb}
\usepackage{amsmath}
\usepackage{amsmath,amscd}
\usepackage{amsthm}
\usepackage{enumerate}
\usepackage{cite}
\usepackage{xcolor}

\usepackage{youngtab}

\usepackage[section]{algorithm}

\makeindex

\title[A footnote to a footnote]
{A footnote to a footnote to a paper of B. Segre}
\date{}

\newcommand{\C}{\mathbb{C}}
\newcommand{\Z}{\mathbb{Z}}
\newcommand{\Pj}{\mathbb{P}}

\newcommand{\T}{\mathbb{T}}

\newcommand{\sym}{\mathrm{Sym}}

\newcommand{\Oc}{\mathcal{O}}

\newcommand{\M}{\mathcal{M}}
\newcommand{\Rm}{\mathcal{R}}

\newcommand{\rank}{\operatorname{rank}}

\newcommand{\pd}[1]{\partial_{#1}}
\newcommand{\Cat}{\mathcal {C}}
\newcommand{\lsp}[1]{\langle {#1} \rangle}

\def\rig#1{\smash{ \mathop{\,\longrightarrow\,}\limits^{#1}}}

\newtheorem{theorem}{Theorem}[section]
\newtheorem{proposition}[theorem]{Proposition}
\newtheorem{lemma}[theorem]{Lemma}
\newtheorem{corollary}[theorem]{Corollary}

\theoremstyle{definition}
\newtheorem{definition}[theorem]{Definition}

\newtheorem{example}[theorem]{Example}
\newtheorem{remark}[theorem]{Remark}


\author[L.~Chiantini]{Luca Chiantini}
\address{Luca Chiantini. Dipartimento di Ingegneria dell'Informazione e Scienze Matematiche, Universit\`a di Siena, Italy}
\email{luca.chiantini@unisi.it}
\author[G.~Ottaviani]{Giorgio Ottaviani}
\address{Giorgio Ottaviani. Dipartimento di  Matematica e Informatica `Ulisse Dini', Universit\`a di Firenze, Italy}
\email{giorgio.ottaviani@unifi.it}
\thanks{}

\begin{document}

%\begin{abstract}

%\end{abstract}

\maketitle

Dedicated to Ciro Ciliberto, for his 70th birthday. \hfill
\medskip

\hfill {\it Si canimus silvas, silvae sint consule dignae.}

\hfill {\it (Vergilius)}
\bigskip

\section{Introduction}\label{sec:intro}

The paper is devoted to a detailed study of sextics in three variables having a decomposition as a sum of nine powers of linear forms. This is
indeed the unique case of a Veronese image $X$ of the plane which, in the terminology introduced by Ciliberto and the first author in 
\cite{CCi02a}, is {\it weakly defective}, and non-identifiable: a general sextic of the $9$-secant variety of $X$ has two minimal decompositions. 

The title originates from a famous paper of 1981,  \cite{ArbarelloCornalba81},  where 
Arbarello and Cornalba state and prove a result on plane curves with preassigned singularities, which is relevant to extend 
the studies of B. Segre on special linear series on curves. The result (Theorem  3.2) says that the linear system of sextics with $9$ 
general nodes in $\Pj^2$ is the unique 
non-superabundant system of plane curves  with general nodes whose (unique) member is non-reduced.

The result is of course relevant also in the theory of interpolation, and in the study of secant varieties to Veronese varieties, with consequences for the
Waring decomposition of forms. The matter turns out to be strictly related to
 the uniqueness of a minimal Waring decomposition of a general form. From this point of view  Theorem  3.2 of \cite{ArbarelloCornalba81}
and Theorem 2.9 of \cite{CCi06} imply that the unique case of general ternary forms of fixed (Waring) rank with a finite number greater than one of
minimal decomposition holds for degree $6$ and rank $9$. Since the unique sextic with 9 general nodes is twice an
elliptic curve, it turns out by proposition 5.2 of  \cite{CCi06} that a general ternary sextic of rank $9$ has exactly two minimal decompositions.
A complete list of cases, in any number of variables, in which the previous phenomenon appears is contained in \cite{COttVan17a}.

Remaining in the case of ternary sextics, already in \cite{ArbarelloCornalba81} the authors observe that if the $9$ nodes are in special position (e.g. when
they are complete intersection of two cubics), then the linear system of nodal curves has also reduced members.
Our analysis starts here, and aims to describe the decompositions of specific ternary sextics, with respect to the postulation
of the corresponding sets of projective points in $\Pj^2$.

We consider a fixed sextic form $F$ in the polynomial ring  $R=\C[x_0,x_1,x_2]$ in three variables. In order to 
effectively produce decompositions of $F$, the first natural step is to consider the {\it apolar} ideal $F^\perp$ of $F$. 

In general, the apolar ideal $F^\perp$ of a form $F$ of degree $d$ is defined via the natural action
of the dual ring $R^\vee$ on $R$. $F^\perp$ is the ideal of elements $D\in R^\vee$ that kill $F$. A classical theorem by Sylvester says that a set of
points is a decompositions of $F$ if and only if the corresponding ideal in $R^\vee$ is contained in $F^\perp$. Indeed, in order to find properties of the 
decompositions of a specific form $F$, it is often sufficient to look at 
homogeneous pieces of  the ideal $F^\perp$, which correspond to the kernels of the (catalecticant) maps  $(R^\vee)_k\to R_{d-k}$ induced by $F$,
and their projective versions $\Cat^k_F$. In \cite{IK}, as well as in \cite{OedOtt13}, it is explained
that if  the image of $\Cat^k_F$ has the expected dimension and cuts the corresponding Veronese variety in $\Pj(R_{d-k})$ in a finite set $A$,
then the set $A$ determines a decomposition of $F$, and many minimal decompositions can be found in this way.

Consider, in particular, the case of a sextic $F$ of rank $8$. In this situation, the image of the catalecticant map of order $3$
determines a linear space of codimension $2$ in $\Pj(R_{3})$, which cuts the $3$-Veronese of $\Pj^2$ in $9$ points. The set $Z$ of $9$ points
is the image of a complete intersection of two cubics in $\Pj^2$, via the Veronese map. Then the $8$ points of a minimal 
decomposition of $F$ are among the points of $Z$.
It is quite easy to find the coefficients of $F$ with respect to the points of $Z$, and then determine the $8$ points that give a minimal decomposition.
The procedure shows that there are sextic forms which lie in between forms of rank $8$ and forms of rank $9$. These are
forms of rank $9$, with a minimal decomposition $Z$ coming from a complete intersection of two plane cubics. Forms of this sort fill
a subvariety $W$ of the secant variety $S^9(X)$ ($X$ being the $6$-Veronese image of $\Pj^2$) which contains $S^8(X)$.

We show that a general sextic in $W$ has a unique decomposition with $9$ powers of linear forms.
From a certain point of view, this is quite surprising. A general sextic $F'$ of rank $9$ has two minimal decompositions of length $9$, coming from
the existence of an elliptic normal curve in $X$ which spans $F'$. When $F'$ sits in $W$, then there is a pencil of elliptic curves
in $X$ which span $F'$, so one may expect infinitely many minimal decompositions, one for each elliptic curve.

The phenomenon can be investigated in terms of the map  $s^9$ from the abstract secant variety $AS^9(X)$ to $S^9(X)$, which is
generically 2:1.
The component of the ramification locus $\Rm$ is strictly connected with the theory of  Terracini loci, introduced in \cite{BallC}. 
Roughly speaking, Terracini loci contain finite
  subsets of $X$ such that the tangent spaces to the points are not independent. A subset $A\subset X$ consisting of $9$ points belongs to 
  the Terracini locus $\T_9(X)$ exactly when the linear span of $v_6(A)$ is contained in $\Rm$. Thus,  when a sextic
  form $F$ has a decomposition $A$ complete intersection of two cubics, which obviously lie
  in the Terracini locus, then in a neighborhood of $F$ the secant  map has degree $2$,
  so it has degree $2$ also in $F$, but the fiber is non-reduced. 
    
  The subvariety $W$ of $S^9(X)$ lies in between $S^9(X)$ and $S^8(X)$, and it is easy to detect because it is easy to compute the rank and the kernel
  of the catalecticant map $\Cat^3_F$. There are other relevant varieties of $S^9(X)$, which can be studied in terms of the geometry of the
  decomposition. Given a general set $A$ of $9$ points in $\Pj^2$, and a general sextic $F$ in the span of $v_6(A)$, a second decomposition of $F$ comes out by
taking a set of $9$ points $B$ linked to $A$ by a cubic and a sextic. In this case, as explained in \cite{AngeC20}, the sum $I_A+I_B$ 
of the homogeneous ideals of $A$ and $B$ determines a subspace of $R_6$, whose orthogonal direction gives coefficients for the (unique)
form $F$ in the intersection $\lsp{v_6(A)}\cap \lsp{v_6(B)}$. In general, the cubic is uniquely determined by $A$ but the sextic moves and then
the set of forms  defined by $ (I_A+I_B)_6$ dominates $\lsp{v_6(A)}$. On the contrary, when $A$ is complete intersection, then also $B$
is complete intersection. We obtain that starting from a complete intersection $A$ the linked sets $B$ determine only a subvariety of forms in
$\lsp{v_6(A)}$, each having a $2$-dimensional set of decompositions. When $A$ varies in the Hilbert stratum of complete intersections of
type $(3,3)$ in $\Pj^2$, we obtain a subvariety $W'\subset W$ of sextics in $S^9(X)$ with a $2$-dimensional set of decomposition, which can be easily described 
because for such forms $F$  the kernel of $\Cat_F^3$ has dimension $3$, so that $W'$ does not contain $S^8(X)$, but
properly contains $S^7(X)$.

Let us notice that the picture is completed by another subvariety $\Rm'$ of the ramification locus, given by forms $F$ with two different and linked
decompositions $A,B$, both sets of nodes of irreducible sextics. The variety $\Rm'$ parametrizes forms with a $1$-dimensional family
of decompositions, and its intersection with $W$ contains $W'$. The variety $\Rm'$ is the most mysterious object in the picture, for it cannot
be completely described in terms of the catalecticant map $\Cat^3_F$ and its kernel.
While the third section is devoted to the geometry (and, somehow, parametric) description of loci in $S^9(X)$, in the last section
we provide for most of them a set of equations. Equations for $S^9(X)$, $W$, $W'$, can be easily described generically by the
vanishing of appropriate minors of the catalecticant map $\Cat^3_F$. 

Equations for $S^8(X)$ and $S^7(X)$ require the determinant of a $27\times 27$
 matrix $A_f$  which represents a flattening, introduced in \cite{LandOtt13}. 
 It is remarkable that the invariant $H_{27}=\det A_f$, joint with the catalecticant matrix
 $\Cat^3$, allows to describe the equations of all the secant varieties to $X=v_6(\Pj^2)$. The ramification locus turns out to be described also
 in terms of $\det(A_f)$ (see Theorem 4.4). We wonder if a description of equations for $\Rm'$ can be achieved in terms of minors of $A_f$. 

\section{Preliminaries}\label{sec:notation}

\subsection{Notation}
In this section, as in most of the paper, we write $\Pj^n$ for the projective space of linear forms in $n+1$ variables, with
complex coefficients.

{\it In the sequel, by abuse, we will identify a form $F$ with the point in the projective space $\Pj(\sym^d(\C^{n+1}))$ associated to $F$,
and also with the hypersurface of equation $F=0$.} 

Fix $N=-1+\binom {n+d}n$ and identify $\Pj^N$ with  $\Pj(\sym^d(\C^{n+1}))$, 
the projective space of degree $d$ forms in $n+1$ variables.

If we denote with $R$ the polynomial ring $R=\C[x_0,\dots, x_n]$, then $\Pj^N$ is  the projective space over the
degree $d$ piece $R_d$, while $\Pj^n$ is  the projective space over  $R_1$.

We denote with $v_d: \Pj^n\to \Pj^N$ the $d$-th Veronese map of $\Pj^n$ which, in the previous notation, maps $L\in \Pj^n$ to $L^d$.
\smallskip

A (Waring) expression of length $r$ of $F$ is an equality
$$  F = a_1L_1^d+ \dots + a_rL_r^d$$
where the $L_i's$ are linear forms and the $a_i$'s are complex coefficients. Each $L_i$ represents a point in $\Pj^n$.
The expression is {\it non-redundant} if the $L_i^d$ are linearly independent in $\sym^d(\C^{n+1})$, and all the coefficients are non-zero.
Non-redundant means that one cannot find a proper  subset $A'\subset A$ such that also $A'$ is a decomposition of $F$. 

We are aware that, since we are working over the algebraically closed field $\C$, we could get rid of the coefficients $a_i$'s
in a Waring expression of $F$. It is convenient for us to maintain the coefficients, because we will compare,
below, Waring expressions in which the linear forms $L_i$'s are fixed and the coefficients move.

Given a Waring expression of $F$, we call (Waring) {\it decomposition} of $L$ the finite set $A=\{L_1,\dots,L_r\}\subset\Pj^n$.
The length of the decomposition $A$ is the cardinality of the finite set $A$, often denoted with $\ell(A)$.\\ 
From the geometric point of view, $A$ is a decomposition of $F$
if and only if $F$ belongs to the linear span $\lsp{v_d(A)}$ of $v_d(A)$, and $A$ is non-redundant if one cannot find a proper  
subset $A'\subset A$ such that also $A'$ is a decomposition of $F$.

The {\it rank} of $F$ is the minimal $r$ for which $F$ has a Waring expression of length $r$. The {\it border rank} of $F$ is the minimal
$r'$ such that $F$ is limit of forms of rank $r'$.

\subsection{Catalecticant maps}

The {\it projective catalecticant map} associated to a form $F$, as suggested in \cite{IK}, is defined as follows.

Denote with $R^\vee=\C[\pd 0,\dots,\pd n]$ the ring of linear operators. There is a natural contraction map $\Cat : R^\vee\times R \to R$
defined by linearity and  by
$$ \Cat ( \pd {i_1}\cdots\pd {i_k}, F) \mapsto  \frac {\partial ^k F}{\pd {x_{i_1}}\cdots \pd{x_{i_k}}}.$$
Notice that  $\Cat$ maps $(R^\vee)_k \times R_d$ bilinearly to $R_{d-k}$.

For  fixed $F\in R_d$ and $k\leq d$, define the {\it catalecticant map of order $k$} associated to $F$ as
$$ \Cat^k_F: (R^\vee)_k\to R_{d-k}  \qquad  D \mapsto \Cat(D,F).$$

Call {\it $k$-polar space of $F$} the subspace $Im(\Cat^k_F)$ of $R_{d-k}$, and call {\it $k$-polar projective space of $F$}
its projectification $P^k_F=\Pj(Im(\Cat^k_F))$.
\smallskip

We will use in the sequel the correspondence between $R_d$ and $R^\vee_d$ induced by the monomial basis.

Properties of catalecticant maps in connection with decompositions of a form $F$ are deeply studied in \cite{IK}. We recall some basic facts.

\begin{theorem}\label{apolar} If $A=\{L_1,\dots,L_r\}\subset \Pj^n$ is a decomposition of $F$, then the 
projective polar space $P^k_F$ is contained in the linear span of $L_1^{d-k},\dots,L_r^{d-k}$.\\
Thus,  $\dim(P^k_F)\leq \rank(F)$.
\end{theorem} 

Comparing dimensions, one obtains.

\begin{corollary}\label{uguale} In the setting of Theorem \ref{apolar}, if moreover $P^k_F$ has (projective) dimension
$r-1$, then the points $L_1^{d-k},\dots,L_r^{d-k}$ are linearly independent, and 
   $$P^k_F= \lsp {L_1^{d-k},\dots,L_r^{d-k}}.$$
   In particular, if $\dim(P^k_F)=\rank(F)$, then any minimal decomposition $A$ of $F$ satisfies
   $$v_{d-k}(A)\subset P^k_F\cap v_{d-k}(\Pj^n).$$
\end{corollary}

We notice that it is easy to construct the matrix of $\Cat^k_F$ with respect to the monomial basis, and then check the dimension 
of the image of $\Cat^k_F$. 

The Apolarity Theorem (see \cite{IK} Section 4) determines a fundamental link between catalecticant maps and decompositions.
We will use it in the form below

\begin{theorem}\label{apolarem} 
The direct sum of the kernels of the catalecticant maps $\Cat^k_F$ is an (artinian, homogeneous) ideal in $R^\vee$, called
the \emph{apolar ideal} $F^\perp$ of $F$. If the homogeneous ideal $I_{A^\vee}$ of a finite set $A^\vee$, in the projective space $(\Pj^n)^\vee$ 
associated to $R^\vee$, is contained in $F^\perp$, then the corresponding set $A\subset\Pj^n$ determines a decomposition of $F$.\\
\end{theorem}

Notice that when $I_{A^\vee}\subset  F^\perp$, then the $k$-polar space of $F$ is contained in the span of $v_{d-k}(\Pj^n)$ for all $k$.
It is well known that the converse fails. Thus in general one cannot hope to find a decompositions of $F$ by looking only to
one catalecticant map. Yet, we will see that there are cases in which just one catalecticant map is sufficient.

\begin{remark}\label{kerim} 
Consider an element $D$ of the kernel of  $\Cat^k_F$, and an element $D'$ of the kernel of  $\Cat^{k'}_F$.
We can consider $D,D'$ as polynomials in $R^\vee$, with degrees $k,k'$ resp. Assume that $D,D'$ have no common factors.\\
 %Let $G,G'\in R$ be the polynomials  corresponding to $D,D'$ (in the duality given by the monomial basis).
  By the Apolarity Theorem \ref{apolarem} the complete intersection
 of $D,D'$ determines a decomposition $Z$ of $F$, of length $kk'$.\\
Since $F^\perp$ is strictly bigger than the ideal generated by $D,D'$, for  $F^\perp$ is artinian, one can have decompositions $Z'$ of $F$ which are different,
even disjoint, from $Z$.
\end{remark}

\subsection{Hilbert functions}

Let $Z$ be a finite, reduced subset of length $r$ in $\Pj^n$.  Fix a set  of representatives for the coordinates points of $Z$. 
The evaluation of polynomials at the chose representatives provides for any degree $d$
a linear map $\rho_r:R_d\to  \C^r$ whose image has dimension which depends only on $Z$ and not on the choice of the representatives.
The Hilbert function of $Z$ is the map
$$ h_Z:\Z\to \Z\qquad h_Z(d)=\dim Im(\rho_d).$$
We will often use the difference $Dh_Z(d)=h_Z(d)-h_Z(d-1)$, and the h-vector of $Z$, which is the $t$-uple of non-zero values of $Dh_Z$.
\smallskip

Clearly $h_Z(i)=Dh_Z(i)=0$ when $i$ is negative. Also   $h_Z(0)=Dh_Z(0)=1$. Moreover $h_Z(d)$ is non decreasing, and $h_Z(d)=r$ for all large $d$.
Consequently $Dh_Z$ is non-negative and $Dh_Z(d)=0$ if $d$ is sufficiently large.\\
 Other standard properties of $h_Z$ and $Dh_Z$ are known  in the literature, and are collected  e.g. in \cite{C19}. We list the most relevant of them.
\begin{itemize}
\item[a1)] for all $d\gg 0$, $\sum_{i=0}^d Dh_Z(i) = h_Z(d)=r$;
\item[a2)] if $Dh_Z(i)=0$ for some $i>0$, then $Dh_Z(j)=0$ for all $j\geq i$;
\item[a3)] if $Z'\subset Z$ then for all $i$,  $h_{Z'}(i)\leq h_Z(i)$ and $Dh_{Z'}(i)\leq Dh_Z(i)$. 
\item[a4)] if $Z$ is contained in a plane curve of degree $q$, then $Dh_Z(i)\leq q$ for all $i$. 
\end{itemize}

The application of Hilbert functions to tensor analysis is based on the trivial observation that $h_Z(1)<\min\{n+1,r\}$, i.e. the evaluation map
in degree $1$ is not of maximal rank, if and only if $Z$ is linearly dependent and fails to span $\Pj^n$. This implies that when
$\binom{n+d}n\geq r$, then $h_Z(d)<r+1$ if and only if $v_d(Z)$ is linearly dependent.\\
It follows that if $A,B$ are two different, non redundant decompositions of the same form $F\in R_d$, and we take $Z$ to be the union $Z=A\cup B$,
then $h_Z(d)<r$, so that 
\begin{itemize}
\item[b0)] $Dh_Z(d+1)\geq 0$. 
\end{itemize}

The next properties of the union $Z=A\cup B$, which will be used in the sequel, follow from results in
the theory of finite subsets in projective space, as the Cayley-Bacharach theorem, or the Macaulay Maximal Growth principle. 
For an account, we refer to \cite{C19} or \cite{AngeC20}.

\begin{proposition}\label{hf} Let $A,B$ be two different, non-redundant decomposition of a form $F$of degree $d$ in $n+1$ variables. Put $Z=A\cup B$
and assume $A\cap B=\emptyset$.
\begin{itemize}
\item[b1)]  For all $j\leq d+1$ one has
$$ \sum_{i=0}^j Dh_Z(i) \leq \sum_{i=d+1-j}^{d+1} Dh_Z(i).$$
\item[b2)]  The projective dimension of the intersection $\lsp{v_d(A)}\cap \lsp{v_d(B)}$ is equal to $(\sum_{i>d} Dh_Z(i))-1$. In particular,
if $Dh_Z(d+1)=1$ and $Dh_Z(d+2)=0$, then the intersection contains only $F$.
\item[b3)]  Assume $n=2$, i.e. $F$ is a ternary form. If for some $0<i<j$ one has $Dh_A(i)>Dh_Z(j)>0$ then $Dh_Z(j+1)<Dh_Z(j)$.
\end{itemize}
\end{proposition}

\subsection{Ramification and the Terracini locus}

Following \cite{BallC}, Terracini loci can be defined as follows.

\begin{definition}\label{Ter}  The Terracini locus $\T_r(X)$ of a variety $X\subset \Pj^N$ is the closure of the locus in $X^{(r)}$ of subsets $A$
of cardinality $r$ of the regular part $X_{reg}$, such that the span of the tangent spaces to $X$ at the points of $A$ has dimension
smaller than the expected value $r(n+1)-1$.
\end{definition}

The locus $\T_r(X)$ is connected to the structure of the natural map (projection)  from the abstract secant variety $AS^r(X)$ and $S^r(X)$,
$$ s^r:AS^r(X)\to S^r(X) .$$\\
Given a linearly independent finite set $Y\subset X_{reg}$, by Terracini's Lemma and its proof the span of the tangent 
spaces to $X$ at the points of $Y$ correspond to the image of the differential of $s^r$. Since for a general $P\in\lsp Y$,  $ AS^r(X)$
is smooth of dimension $r(n+1)-1$ at $(P,Y)\in  AS^r(X)$, then {\it
$Y$ belongs to $\T_r(X)$ if and only if the linear span $\lsp Y$ is contained in the ramification locus of the map $s^r$.}

\section{Plane sextics}\label{sec:6}

In this section we will analyze in details the situation of forms of degree $6$ in three variables, i.e. sextic plane curves. 
We will see how the study of the catalecticant map determines  several loci in the secant varieties of the surface $X=v_6(\Pj^2)$.

\begin{remark}\label{in}  Sextics in three variables are parameterized by $\Pj(\sym^6(\C^3))=\Pj^{27}$.\\
  By \cite{AlexHir95}, the (Waring) rank of a general ternary sextic is $10$. \\
  Sextics of (border) rank $9$ determine a hypersurface
in the space of sextics:  the $9$-secant variety to the $6$-Veronese of $\Pj^2$.
Moreover, by \cite{ArbarelloCornalba81} and by \cite{CCi06}, a general form of rank $9$ has exactly two minimal decompositions.\\
Sextics of rank $8$ determine an irreducible subvariety of dimension $23$ in $\Pj^{27}$, i.e. the $8$-secant variety to $v_6(\Pj^2)$. By 
\cite{COttVan17a}, a general sextic of rank $8$ has a unique decomposition of length $8$.
\end{remark}

Following \cite{OedOtt13}, let us look what happens to the catalecticant map of order three applied to a sextic $F$ of rank $<10$.
We will focus mainly on the cases in which $F$ has rank either $8$ or $9$.
\smallskip

We start by recalling the following, elementary fact.

\begin{remark} Let $C$ be a form in the kernel of $\Cat_F^k$. Then $C$ corresponds by duality to a form in $R$ (that we continue to denote as $C$),
and its coefficients correspond to the orthogonal to a hyperplane in $\Pj(R^\vee_k)$ that contains $P^k_F$. Since
the Veronese map  and the duality are both given in terms of the monomial basis, then the intersection  of $P^k_F$ with $v_k(\Pj^2)$
corresponds to a divisor that defines $v_k$, i.e. to a curve of degree $k$ in $\Pj^2$, which is exactly  $C$. Thus $P^k_F \cap v_k(\Pj^2)= v_k(C).$
\end{remark}

\begin{proposition}\label{rank9}  If $F$ has rank $r<10$ then $\dim(P^3_F)\leq r-1$. If $\dim(P^3_F)= 8$ and the kernel of $\Cat^3_F$
is generated by an irreducible form $C$, then $F$ has border rank at most $9$, and rank $\geq 9$.
\end{proposition}
\begin{proof}
 If we take a minimal decomposition  $A$ of $F$, then the ideal of $A$ certainly contains $10-r$ independent cubics. Thus
 from Theorem \ref{apolar} we know that the kernel of $\Cat^3_F$ has dimension  at least $10-r$, so that $\dim(P^3_F)\leq r-1$.\\
Now assume that  the dimension of $P^3_F$ is exactly equal to $8$, so that $F$ cannot have rank smaller than $9$.
   $P^3_F$ is then a hyperplane in  $\Pj(\sym^3(\C^3))=\Pj^9$ which intersects $v_3(\Pj^2)$ in $v_3(C)$. 
   Since $F^\perp $ certainly contains some form $C'$
   independent from $C$, and by the Apolarity Theorem  $C\cap C'$ is a decomposition of $F$, hence $F$ lies in the span
   of $v_6(C)$, which is the Veronese image of an irreducible plane cubic of arithmetic rank $1$. Since irreducible
   curves are never defective, then $F$ has border rank at most $9$ with respect to $v_6(C)$, hence also with respect to $v_6(\Pj^2)$.
\end{proof}

 Notice that \emph{not every} general set of $9$ points in $C$ determines a decomposition of $F$. So, it is not sufficient that the span
 of $v_3(A)$ contains $P^3_F$ for $A$ to be a decomposition of $F$.\\

\begin{example}\label{6rank8}
Let $F$ be a sextic of rank $8$. Since $8$ points always lie in a pencil of cubic curves, then $\dim(P^3_F)\leq 7$. 
Assume that $P^3_F$ has dimension exactly $7$
and  look at the intersection of $P^3_F$ with the $3$-Veronese surface $S=v_3(\Pj^2)$.\\
Since $\deg(S)=9$, it is clear that $P^3_F$ intersects $S$ either in a curve $\Gamma$, or in a scheme of length $9$.\\
The former case cannot happen. Namely the curve $\Gamma$ would be the image in $v_3$ of a plane curve whose ideal in degree $3$
has linear dimension $1$, but no such plane curves exist.\\
Thus  $Z=P^3_F\cap S$ is a subscheme of length $9$. It is easy to see (\cite[Section 4]{OedOtt13}) that for a general choice of the form 
$F$ the scheme $Z$ is reduced, hence it consists of $9$ distinct points. \\
{\it It follows that any  decomposition of length $8$ of $F$ is  contained in $Z$}, for any set of length $8$ sits in two independent cubics,
and by the Apolarity Theorem  these two cubics are in the kernel of $\Cat^3_F$. \\
When $Z$ is reduced, it consists of a minimal decomposition, plus one extra point.
\end{example}

\begin{proposition}\label{2cub2} Let $F$ be a sextic ternary form such that $\dim(P^3_F)=7$ and $Y=P^3_F\cap S$ is reduced. Assume that $F$
has rank $8$. Then there exists only one decomposition of $F$ of length $8$, i.e. $F$ is identifiable, in the sense of \cite{CCi06}.
\end{proposition} 
\begin{proof} Put $Z=v_3^{-1}(Y)$. By Example \ref{6rank8} every decomposition of length $8$ of $F$ is a subset of $Z$. 
Since $Z$ is complete intersection of two cubics,
then  $Z$ imposes independent conditions to curves of degree $\geq 3+3 -2 =4$.
In particular, $v_6(Z)$ is linearly independent, thus every subset of $v_6(Z)$ is linearly independent.\\
Assume that there are two different decompositions $A,A'$ of length $8$. Since $A,A'\subset Z$ by Remark \ref{6rank8}, then $B=A\cap A'$ has  
length $7$. Since $v_6(A)$ and $v_6(A')$ are linearly independent, then the spaces $\lsp{v_6(A)}$ and $\lsp{v_6(A')}$, both containing $F$,
meet exactly in $\lsp{v_6(B)}$. Then $F\in\lsp{v_6(B)}$, which means that $\rank(F)\leq 7$, a contradiction.
\end{proof}

In conclusion, we get
\begin{corollary}\label{2cub1} Let $F$ be a sextic ternary form such that $\dim(P^3_F)=7$ and $Y=P^3_F\cap S$ is reduced. Then 
$Z=v_3^{-1}(Y)\subset\Pj^2$ is a complete intersection of two cubics, and either
\begin{itemize}
\item $F$ has rank $8$ and a unique decomposition of length $8$, contained in $Z$; or
\item $F$ has rank $9$, and $Z$ is a minimal decomposition of $F$, i.e.
there exists a  minimal decomposition of $F$ which is a complete intersection of $2$ cubics.
\end{itemize}
\end{corollary}

An algorithm that can guarantee that a ternary form $F$ of degree $6$ has rank $8$, and find a decomposition, is the following (see  \cite{OedOtt13}).\\
Compute $\dim(P^3_F)$. If $\dim(P^3_F)>7$, then $F$ cannot have rank $8$. Assume that $\dim(P^3_F)=7$ (this will be true for general forms of rank $8$).\\
Compute two generators $D,D'$ of the kernel of $\Cat^3_F$ and compute the coordinates of the points of the complete
intersection $Z$ of $D$ and $D'$ (this can require approximation). Assume that $Z$ has length $9$ (this will be true for general forms of rank $8$).\\
Then $Z=\{L_1,\dots,L_9\}$. If we identify each $L_i$ with a linear form in $R$, then there exists a linear combination
$F=a_1L_1^6+\dots + a_9L_9^6$ and we know that the linear combination is unique, for $v_6(Z)$ is linearly independent. Compute
the  $a_i$'s.\\
By Proposition \ref{2cub2}, $F$ has rank $8$ if and only if one coefficient $a_i$ is $0$. In this case, by dropping the summand $a_iL_i^6$,
we also obtain the unique decomposition of length $8$ of $F$.
\medskip

It is known that since  a general point in $\Pj^{17}$ has exactly  two different decompositions with respect to  an elliptic normal curve,
then a general sextic of rank $9$ has exactly two minimal general decomposition. \\
A geometric way to find the second decomposition, once one of them is known, can be described in terms of liaison.\\
We need first a series of preparatory results on the interaction between Hilbert functions and decompositions.

\begin{lemma}\label{IL} ({\bf Intersection Lemma}) Let $F\in R_d$ be a form with two non-redun-dant decompositions $A,B$, of length $\ell(A)\geq \ell(B)$. Assume 
$A\cap B\neq \emptyset$. Then there exists a form $F'\in R_d$ with two non-redundant decompositions $A'\subset A$, $B'\subset B$ such that $A'\cap B'=\emptyset$
and $\ell(A')\geq \ell(B')$, and moreover either $\ell(A')> \ell(B')$ or $A'\neq A$.
\end{lemma}
\begin{proof} Let $A=\{L_1,\dots,L_r\}$ and $B=\{L_1,\dots,L_j,M_{j+1},\dots, M_{r'}\}$ with $L'_i\notin A$ for $i=j+1,\dots,r'$ ($r=\ell(A), r'=\ell(B)$). Thus
$A\cap B=\{L_1,\dots,L_j\}$ and $j>0$. Write
\begin{gather*} F= a_1L_1^d+\dots +a_jL^d_j +a_{j+1}L^d_{j+1}+\dots+ a_rL^d_r \\
	F= b_1L^d_1+\dots +b_jL^d_j +b_{j+1}M^d_{j+1}+\dots+ b_{r'}M^d_{r'}.
\end{gather*}
Since $A,B$ are non-redundant, then $v_d(A),v_d(B)$ are both linearly independent, and the $a_i$'s and $b_i$'s are all non-zero. Consider the
form $F'$ obtained by subtracting from $F$ the first $j$ summands of the second Waring expression above
  \begin{gather*} F'= (a_1-b_1)L_1^d+\dots +(a_j-b_j)L^d_j +a_{j+1}L^d_{j+1}+\dots+ a_rL^d_r \\
	F'= \phantom{(a_1-b_1)L_1^d+\dots +(a_j-b_j)L^d_j +}b_{j+1}M^d_{j+1}+\dots+ b_{r'}M^d_{r'}.
\end{gather*}
Then $B'=\{M_{j+1},\dots,M_{r'}\}$ is a decomposition of $F'$, of length $r'-j<r'$, which is non-redundant since $b_{j+1},\dots,b_{r'}\neq 0$ and $v_d(B')$
is linearly independent. Let $A'=A\setminus\{L_i: i\leq j$  and $ a_i=b_i\}$. Then also $A'$ is a decomposition of $F'$, which is non-redundant since 
$v_d(A')$ is linearly independent and $a_{j+1},\dots,a_r\neq 0$. $A'$ and $B'$ are clearly disjoint, and $\ell(A')\geq r-j\geq r'-j=\ell(B')$. \\
Finally, if $A'=A$ then $\ell(A')=r>r'-j$.
\end{proof}

The next result tells us where we must look, geometrically, in order to find minimal decompositions of a sextic form.\\
We give a long proof in which we analyze one by one the several cases, just because we want to stress the use of the Hilbert
function for tensor analysis. Indeed, some steps could be shortened by using \cite{Ball19} or Section 3 of \cite{AngeC20}.
See also \cite{MourOneto20} for an approach to the problem.

\begin{proposition}\label{rleq9} Let $F$ be  a ternary sextic with two different non-redundant decompositions $A,B$, with length $\ell(A)\leq 9$
and $\ell(B)\leq \ell(A)$. Assume that no $4$ points of $A$ are aligned, and that $A$ does not lie in a conic.\\ Then $\ell(B)=\ell(A)=9$ and $A\cap B=\emptyset$.\\
Furthermore $A\cup B$ lies in a cubic curve $C$, and it is complete intersection of $C$ and a sextic curve.
\end{proposition}
\begin{proof} We set $r=\ell(A), r'=\ell(B)$ and we make induction on $r,r'$. Put $Z=A\cup B$ and notice that $Dh_Z(7)>0$, by b0). \\
 If $r\leq 3$, for any $r'$ we have a contradiction with \cite{BucGinenskyLand13} Corollary 2.2.1. \\
 First assume that $A,B$ are disjoint.\\
Let $r=4,5$. Since four points of $A$ are not aligned, In both cases   $\sum_{i=0}^2 Dh_Z(i)\geq r$,
hence $\sum_{i=5}^7 Dh_Z(i)\geq r$ by b1). It follows that $Dh_Z(3)=0$, which contradicts $Dh_Z(7)>0$, by a2).\\
Assume $r=6$,  so that $\ell(Z)\leq 12$. Since $A$ is not aligned, then $Dh_A(1)=2$. If $Dh_A(2)\leq 1$, then A sits in two independent conics, which, by Bezout, 
is possible only if $A$ has $5$ aligned points, a contradiction. Thus $Dh_A(2)\geq 2$, so that the h-vector of $A$ is either $(1,2,3)$, or $(1,2,2,1)$. 
In the former case, arguing as in cases $r=4,5$, we get   $Dh_Z(3)=0$, a contradiction. In the latter case, we have $\sum_{i=0}^3 Dh_Z(i)\geq 6$, thus by b1)
also $\sum_{i=4}^7 Dh_Z(i)\geq 6$. Since $\ell(Z)\leq 12$, we must have $\sum_{i=0}^3 Dh_Z(i)= 6$ which, by a3), implies $Dh_Z(i)=Dh_A(i)$ for $i=0,\dots,3$.
Then $Dh_Z(3)= 1$, so that, by b3), $Dh_Z(4)=0$, which is not consistent with b0).\\
Assume $r=7,8$. $A$ cannot be contained in conics, so that $Dh_A(2)=3$. If $Dh_Z(3)\leq 2$, then by applying b3) for $i=4,5$ we
see that $Dh_Z(5)=0$, which contradicts b0). Thus $Dh_Z(3)\geq 3$, so that $\sum_{i=0}^3 Dh_Z(i)\geq 9$. Then by b1) also $\sum_{i=4}^7 Dh_Z(i)\geq 9$,
which is not consistent with $\ell(Z)\leq 2r$.\\
So, we are left with the case  $r=9$. As above $Dh_Z(2)=3$, and $Dh_Z(3)\leq 2$ yields a contradiction. Thus the h-vector of $Z$ 
starts with $(1,2,3,q,\dots)$ with $q=Dh_Z(3)\geq 3$. Then $\sum_{i=0}^3 Dh_Z(i)= 6+q$, so that by b1) $\sum_{i=4}^7 Dh_Z(i)\geq 6+q$. Since
$\ell(Z)\leq 18$, it follows $q=3$. This means that $\ell(Z)=18$, i.e. $\ell(B)=9$, and moreover $h_Z(3)=9$, so $Z$ lies in a cubic curve.\\
By  Proposition  \ref{hf}, one computes that the unique possibility for the h-vector of $Z$ is $(1,2,3,3,3,3,2,1)$. 
By \cite{AngeCMazzon19} Lemma 5.3 we know that $Z$ has the Cayley-Bacharach property. Since the h-vector of $Z$ is the same than the
h-vector of a complete intersection of type $(3,6)$, by \cite{Davis84}  $Z$ is complete intersection of a cubic and a sextic curve.
Then $A$ and $B$ are linked by a cubic and a sextic curve.
Finally, the case where $A\cap B\neq \emptyset$ implies, by Lemma \ref{IL}, the existence of a sextic $F'$ with  non-redundant decompositions $A'\subset A$
and $B'$, such $\ell(B')\leq \ell(A')$, and either $\ell(A')\leq 8$ or $\ell(B') \leq 8$. The existence of $F'$ is excluded by the previous argument. 
\end{proof}

 Given a general sextic $F$ of rank $9$, and a decomposition $A$ of $F$, we can find a second decomposition $B$ of
$F$ with a procedure introduced e.g. in \cite{AngeC20}, Section 4.

\begin{remark}\label{second}
Let $F$ be a sextic of rank $9$, with a decomposition $A$ of length $9$.\\
From the kernel of $\Cat^3_F$, which has dimension $1$ is $F$ is general, we find  a cubic $C$ containing $A$.
A resolution of the homogeneous ideal of $A$ is given  by
$$ 0\to R(-5)^3 \stackrel {M}\to R(-4)^3\oplus R(-3) \to I_A\to 0.$$
$M$ is the Hilbert-Burch matrix of $A$, whose minors provide generators for $I_A$. The degrees of $M$ are
$$\begin{pmatrix} 1 & 1 & 1 \\1 & 1 & 1\\1 & 1 & 1\\ 2& 2 & 2 \end{pmatrix}.$$
We know that $A\cup B$ is complete intersection of $C$ with a sextic. The homogeneous ideal of the residue $B$ of $A$ in the intersection
is obtained by erasing the bottom row of $M$ and adding a column of three quadrics. The maximal minors of the matrix  $M'$ that we obtain
generate the homogeneous ideal of a scheme $B$ linked to $A$ in a complete intersection $(3,6)$.\\
The previous procedure produces a set $B$ which decomposes some form $F'\in v_6(A)$. In other words, 
the intersection $\lsp{v_6(A)}\cap\lsp{v_6(B)}$ is a form $F'$ in $v_6(A)$ which depends on the choice of the three quadrics in $M'$.\\
One can find $F'$ by taking the sum $I_A+I_B$, which determines in $R_6$ a linear subspace $H$ of codimension $1$.
The coefficients of $F'$ form a vector orthogonal to $H$.\\
Thus, in order to find the second decomposition of the given $F$, one needs to choose the three quadrics $Q_1,Q_2,Q_3$ in $M'$ so that
the orthogonal to $H$ matches with the coefficient of $F$. In practice,  one needs to solve
a linear system whose matrix has entries linear in the coefficients of the $Q_i$'s. A convenient way to write this linear system
into M2\cite{Macaulay2} is to ask that the minors of $M'$ are apolar to $F$.
\end{remark}

With a procedure similar to the proof of Proposition \ref{rleq9}, we can prove the following property of sextics of rank $9$
with a decomposition $A$ complete intersection of two cubics.

\begin{proposition}\label{2cub3} Let $F$ be a sextic ternary form such that $\dim(P^3_F)=7$ and $Y=P^3_F\cap S$ is reduced. Assume that $F$
has rank $9$. Then $A=v_3^{-1}(Y)$ is the unique decomposition of $F$ of length $9$, i.e. $F$ is identifiable, in the sense of \cite{CCi06}.
\end{proposition} 
\begin{proof} From Proposition \ref{2cub1} we know that $A$ is a complete intersection of two cubics. In particular, no $4$ points of $A$ are aligned
and $A$ does not lie in a conic. The h-vector of $A$ is $(1,2,3,2,1)$.\\
Assume that there exists another decomposition $B$ of length $9$. Then by Proposition \ref{rleq9} we know that $A\cap B=\emptyset$ and
 $Z=A\cup B$ lies in a cubic, so that the h-vector of $Z$ starts with $(1,2,3,3,\dots)$. Since $Z$ lies in a plane curve of degree $3$, 
by a4) $Dh_Z(i)\leq 3$ for all $i$.\\
Next, we claim that $Dh_Z(4)=3$. Indeed if $Dh_Z(4)\leq 2$ then by b3) $Dh_Z(5)\leq 1$ and $Dh_Z(6)=0$, which contradicts b0). If 
$Dh_Z(4)\geq 4$, then $\sum_{i=5}^7 Dh_Z(i)<6$, which contradicts b1). The same argument proves that $Dh_Z(5)=3$, and
$Dh_Z(6)\leq 2$. Then $Dh_Z(6)+Dh_Z(7)=3$. Since by b3) $Dh_Z(6)>Dh_Z(7)\geq 1$, then the h-vector of $Z$ is $(1,2,3,3,3,3,2,1)$.\\
By \cite{AngeCMazzon19} Lemma 5.3 we know that $Z$ has the Cayley-Bacharach property $CB(6)$. Since moreover the h-vector of $Z$ is the same than the
h-vector of a complete intersection of type $(3,6)$, by \cite{Davis84}  $Z$ is complete intersection of a cubic and a sextic curve.
Then $A$ and $B$ are linked by a cubic and a sextic curve. This implies, by \cite{Migliore} 
that the h-vector of $B$ is $(1,2,3,2,1)$. Thus $B$ lies in the base locus of a pencil of cubic curves. Since $B$ is a decomposition of $F$, and $\dim(P^3_F)=7$,
 there is  only one pencil of cubic curves that contain $B$, and the base locus of the pencil is $A$. Comparing the degrees, we get
 $A=B$.
\end{proof}

from a certain point of view, Proposition \ref{2cub3} is rather surprising. Namely, as explained in \cite{ArbarelloCornalba81} Theorem  3.2, 
the $6$-Veronese surface $v_6(\Pj^2)$
is $9$-weakly defective, in the sense of \cite{CCi02a}: a general hyperplane which is tangent to $v_6(\Pj^2)$ at $9$ points is tangent along
an elliptic normal curve $C$. Thus, by \cite{CCi06} Theorem 2,9 and Proposition 5.2,  $F$ has two minimal decompositions of length $9$,
both lying in $C$.
Notice that, in general,  none of the two decompositions  lies in two independent cubic curves.\\
If the form $F$ has a decomposition $Z$ of length $9$ which is a complete intersections of two cubics, then there are infinitely
many elliptic normal curves in $v_6(\Pj^2)$ containing $v_6(Z)$, so one may expect the existence of infinitely many decompositions
of length $9$ of $F$, for any  elliptic normal curves  $v_6(Z)$ could provide a second decomposition $Z'\neq Z$ of $F$.
On the contrary,  in this case the construction collapses and any elliptic normal curves provides only one decomposition for $F$, namely $A$ alone.
\smallskip

A  motivations for the peculiar behavior  of sextics with a decomposition complete intersection of two cubics can be
explained in terms of the ramification of the map $s^9$ from the abstract secant variety $AS^9(v_6(\Pj^2))$ to $S^9(v_6(\Pj^2))$, and the Terracini locus of 
$v_6(\Pj^2)$.

\begin{remark} 
For $X=v_6(\Pj^2)$, by Definition \ref{Ter} a linearly independent set $v_6(A)$ of length $9$ belongs to $\T_9(X)$ if the $9$ tangent planes at the points of 
$v_6(A)$ span a subspace of dimension $\leq 25$. When $A$ is complete intersection of two cubics $\Gamma_1,\Gamma_2$, then certainly 
$v_6(A)$ belongs to the Terracini locus. Indeed, there is a $2$-dimensional
family of sextic curves in $\Pj^2$ which are singular at the points of $A$: the sextics formed by the union of two cubics in the pencil
generated by $\Gamma_1,\Gamma_2$. Any sextic like that corresponds to the intersection of $X$ with a hyperplane which contains the
tangent planes to $X$ at the points of $v_6(A)$. Thus, the span of the tangent planes at the points of $v_6(A)$ has codimension at least $2$
 in $\Pj^N=\Pj^{27}$.\\
By Remark \ref{in}, the map $s^9:AS^9(X)\to\Pj^{27}$ is generically $2:1$, when $A$ is a general set of $9$ points in the plane.
By Proposition \ref{2cub3} we see that general points $(F,v_6(A))$ such that $A$ is complete intersection of two cubics are
in the ramification locus of $s^9$.\\
The last observation is consistent with the fact that the differential of $s^9$ degenerates at points $(F,v_6(A))$ such that $A$ is complete 
intersection of two cubics.
  \end{remark}

  \begin{definition} \label{def:RW}
   We denote with $\Rm$ the closure of the set of sextics with a decomposition formed by a set $A$ of nine points which are nodes
   of an irreducible sextic. By \cite{Harr86},  $\Rm$ is irreducible.\\
   We denote by $W$ the closure of the set of sextics with a decomposition formed by a set $A$ of nine points, complete intersection of two cubics. 
  \end{definition}
  
 \begin{remark} \label{RW} The previous discussion proves that, in the new notation, $W$ properly contains $S^8(X)$, and both $W$ and $\Rm$ are contained in 
 the Terracini locus $\T_9(X)$, which is the closure of the image of the locus of smooth points in  $AS^9(X)$ in which the differential of the 
   map $s^9$ drops rank.  \\
   Indeed, we will show in the next section that $\Rm$ is the unique component of $\T_9(X)$ which intersect
 the smooth locus of $S^9(X)$, while $W$, whose codimension is bigger than $1$,  sits into the singular locus of $S^9(X)$.\\
We will also prove in the next section that $W$ is not contained in $\Rm$, which means that the Terracini locus has at least two components.
\end{remark}

  \begin{remark} For the $6$-Veronese variety $X$ of $\Pj^2$, the Terracini locus $\T_9(X)$ corresponds to sets $A$ of $9$ points which are
  singular in a pencil of sextics.  General sets in $\T_9(X)$ are easy to describe (see Example 5.1 in \cite{BallC}): just take $C$ to be a reduced elliptic
  plane sextic whose singularities are nodes. $C$ has a set $A$ of $9$ nodes which also sit in a cubic curve $D$. All sextics in the pencil 
  generated by $C$ and $D^2$ are singular at $A$. It follows that the $9$ tangent planes to $X$ at the points of $v_6(A)$ span
  a linear space of dimension $25$ in $\Pj^{27}$, the linear space spanned by $X$.\\
  Notice that a $\T_9(X)$ does not contain a general set of $9$ points in the plane, and a general $A\in  \T_9(X)$ is not complete intersection of two cubics.
  \end{remark}
  
General forms contained in the variety $\Rm$ have only one decomposition with $9$ summands.
On the other hand, both $\Rm$ and $W$ contain subvarieties whose elements have many decompositions of length $9$.\\
Namely, if $A$ lies in the Terracini locus $\T_9(X)$, then the liaison procedure introduced in Remark \ref{second} determine a second 
decomposition for some forms in $\lsp{v_6(A)}$.

\begin{lemma}\label{lemprima} Let $C$ be a general (smooth) cubic plane curve and let $A$
be a reduced divisor of degree $9$ on $C$ such that $2A\in|\Oc_C(6)|$, but $A\notin|\Oc_C(3)|$. 
Then there exists a sextic plane curve $G$ which is singular at the points of $A$ and irreducible.
\end{lemma}
\begin{proof} Call $P_1,\dots,P_9$ the points of $A$. For each $i$ call $\epsilon_i$ the
double structure on $P_i$ contained in $C$. For each $P_i$ choose a scheme $\mu_i$ of length $2$
in $\Pj^2$, supported at $P_i$ such that $\mu_i\notin C$. By assumption there exists a sextic $G'$ such that
$G'$ does not contain $C$ and $\epsilon_i\in G'$ for all $i$. Thus, the linear system $\M$ of
sextics in $\Pj^2$ which contain $\{\epsilon_1,\dots,\epsilon_9\}$ has (projective)
dimension $10$. Note that $\M$ contains the double curve $2C$. Since $\epsilon_i\cup \mu_i$ has length $3$,
then every $\mu_i$ imposes at most one condition to curves in $\M$.
Then there exists a sextic $G''$ which contains $\epsilon_i\cup \mu_i$ for all $i$ and $G''$ is different
from $2C$. Since $G''$ contains $\epsilon_i\cup \mu_i$, then $G''$ is singular at the points of $A$.
$G''$  cannot contain $C$, because $A\notin|\Oc_C(3)|$. Thus a general curve $G$ in the pencil generated by
$2C$ and $G''$ is irreducible, by Bertini. 
\end{proof}

\begin{proposition} Let $A$ be a set of $9$ points in $\Pj^2$, which are nodes of a 
general reduced elliptic sextic curve $G$. The generality of $G$ implies that
there exists a unique cubic curve $C$ containing $A$, and moreover $C$ is smooth.
Let $B$ be the residue of $A$ in the intersection of $C$ with a general sextic $S$ 
that passes through $A$. Then also $B$ is the set of nodes of a reduced sextic $G'$.
\end{proposition}
\begin{proof} The first claim is classical. Just to see an argument, observe
that for a general set $D$ of $8$ points in a general cubic curve $C$, the divisor $2D$
on $C$ lies in some divisor $D'$ of the linear series $|\Oc_C(6)|$. $D'-D$ consists of $2$ points, 
that determine a linear series $\mathcal L =g^1_2$. Choose a Weierstrass point $P_9$ of 
$\mathcal L$. Since $\mathcal L$ cannot have $8$ Weierstrass points, then by monodromy for a 
general choice of $D$ we can assume that $P_9\notin D$. Then $2D+2P_9\in |\Oc_C(6)|$ and then
$D\cup \{P_9\}$ is contained in the singular locus of a sextic curve, by Lemma \ref{lemprima}. Since the 
variety that parametrizes sextic plane curves with $9$ singular points is irreducible, 
we see that on a general cubic $C$ we can find a set of $9$ nodes of an irreducible sextic.\\ 
For the second claim, observe that on $C$ the divisor $2A$ belongs to $|\Oc_C(6)|$, and also
the divisor $A+B$ belongs to $|\Oc_C(6)|$. Thus $2A+2B$ belongs to $|\Oc_C(12)|$.
Then also $2B\in|\Oc_C(6)|$, and the claim follows from the generality of $G,C,S$ and
from Lemma \ref{lemprima}.
\end{proof}

In other words, the previous proposition says that  some forms $F$ of rank $9$
with a decomposition $A$ which lies in the Terracini locus $T_9(X)$, have a second decomposition $B$
of length $9$ which also sits in the Terracini locus $T_9(X)$. Thus the subscheme of the abstract
$9$-secant variety which maps to $F$ is supported at two points, but it has length $4$.\\
The closure $\Rm'$ of the locus of forms $F$ as above provides a proper subvariety of $\Rm$,
whose geometry has non-trivial aspects.

\begin{proposition} A general sextic in $\Rm'$ has infinitely many decompositions of length $9$. Thus $\Rm'$ is the locus in $\Rm$
in which the dimension of fibers of the map $s^9:AS^9(X)\to S^9(X)$ jumps.
\end{proposition}
\begin{proof} Fix a general cubic curve $C$ and a general set $A\subset C$ of length $9$, such that there is an irreducible sextic $S$
 singular at the points of $A$. For a general sextic  $S'$ through $A$ the intersection $S'\cap C$ is a set $A\cup B$ of $18$ points, with $B$ disjoint from $A$.
 As explained in  Proposition \ref{rleq9}, the linear spans $\lsp{v_6(A)}$ and $\lsp{v_6(B)}$ meet in a point $F\in\Pj^{27}$ which represents
 a point in $\Rm'$. Since any element of $\Rm'\cap \lsp{v_6(A)}$ arises in this way, we get a surjective map form the projective space $\Pj$
 over $(I_A)_6$ mod multiples of $C$, to  $\Rm'\cap \lsp{v_6(A)}$. Since $A$ is separated by cubics, then  $\Pj$ has dimension  $27-19=8$.
 Since $\lsp{v_6(A)}$ also has dimension $8$, and $\Rm'\cap \lsp{v_6(A)}$ is a proper subvariety, because $\Rm\neq \Rm'$, the claim follows.
\end{proof}

  The same procedure of liaison produces a special subvariety of $W$.
  
  \begin{example}\label{w'} Let $A$ be a set of $9$ points, complete intersection of two cubics $C_1,C_2$ Fix a general cubic curve $C_3$ and consider the set
  $B=C_1\cap C_3$. Then $A\cap B=\emptyset$, and $Z=A\cup B$ is a complete intersection of a cubic and a sextic curve. Thus the h-vector of $Z$ is 
  $(1,2,3,3,3,3,2,1)$, By \cite{AngeC20} Proposition 2.19, it follows that the spans $\lsp{ v_6(A)}$ and $\lsp{v_6(B)}$ meet in exactly one point,
  corresponding to a sextic form $F$.\\
  Since $F$ has a decomposition given by $C_1\cap C_2$, then $F$ corresponds to a point in the $9$-secant variety of $X$ spanned by a set in
  the Terracini locus. Moreover, from the Apolarity Theorem, $C_1,C_2$ lie in the apolar ideal of $F$. For the same reason,
  looking at the decomposition $C_1\cap C_3$ of $F$, one obtains that also $C_3$ lies in $F^\perp$. Thus we have an example of 
  a form $F$ whose apolar ideal contains $3$ independent cubics $C_1,C_2,C_3$.\\
  Notice that, by the generality of $C_3$, the intersection $C_1\cap C_2\cap C_3$ is empty. Thus the image of $\Cat^3_F$, which is a $\Pj^6$, will not cut
  the 3-Veronese surface in $\Pj^9$. In particular, $F$ has not rank $7$.\\
 By the Apolarity Theorem, a decomposition of $F$ can be found by taking any general pair of cubics in the linear system spanned by $C_1,C_2,C_3$. Thus
 $F$ has a $2$-dimensional family of decompositions.
 \end{example}
 
 \begin{definition}
We denote with $W'$ the closure of the locus of sextic  forms $F$ such that the catalecticant map $\Cat^3_F$ has a kernel of dimension $3$.
\end{definition}

A general choice of three cubics $C_1,C_2,C_3$ determines $F\in W'$, by taking the intersection of
 the spans of $v_6(C_1\cap C_2)$ and $v_6(C_1\cap C_3)$ (any general choice of two pairs will be suitable and determine the same $F$).
 Thus  $W'$ has the  dimension of the Grassmannian of planes in $\Pj^9$, i.e. $21$.
 \smallskip

 We stress that there is no way to determine if $F$ belongs to $W'$ only by looking at the  $9$ projective points of a decomposition of $F$.
 Once one  knows that $F$ has a decomposition $A$ which is complete intersection, and one fixes representatives
 $\{L_1,\dots,L_9\}$ for the points of $A$, then the membership of $F$ in $W'$
 depends on the coefficients of a linear combination of powers $L_1^6,\dots,L_9^6$ that determine $F$ in $\lsp{v_6(A)}$.

 Observe that the locus $S^7(X)$ lies in the intersection of $S^8(X)$ with $W'$. We will see in the next section that, at least set-theoretically, 
   $S^7(X)=S^8(X)\cap W'$.

  \begin{example} Many forms $F$ in the subvariety $W'$ defined above can be computed following the procedure introduced in Section 4
  of \cite{AngeC20}, as explained in Remark \ref{second} .\\
Just to see an example of such a form, consider
  $$ F=(x_0x_1x_2)^2  .$$
  One easily sees that $F^\perp$ contains three independent cubics, which are $x_0^3$, $x_1^3$, $x_2^3$.\\
 We have the explicit decomposition of rank $9$ \cite{BucaBucTeit13}
  $$810(x_0x_1x_2)^2=\sum_{p,q=0}^2e^{2\pi i(p+q)/3}\left(x_0+e^{2\pi ip/3}x_1+e^{2\pi iq/3}x_2\right)^6.$$
  In the next section we will compute equations for the $8$-th secant variety $S^8(X)$. By applying these equations to $F$, one realizes that $F$
  does not belong   to $S^8(X)$. Thus $F$ has rank (and border rank) $9$. This proves that $W'$ is not contained in $S^8(X)$. 
  \end{example}
  
  The relations among the subvarieties $\Rm$, $\Rm'$, $W$, $W'$ of $S^9(X)$ reflect the rich geometry of 
   secant varieties, as soon as the genus approaches the generic value.\\
The most complicate object to describe remains $\Rm'$, for which we do not have a set of equations. 
We hope that the geometric structure of $\Rm'$ will be clarified, in a future footnote, maybe.

  \section{Equations for loci in $S^9$.}
   We introduce a Young flattening to get a more refined study of plane sextics.
  Let $T$ be the tangent bundle of $\Pj^2$, weconsider the rank three bundle $E=\sym^2T$.
  The space of sections $H^0(\sym^2T)$ is $27$-dimensional and can be identified with the $SL(3)$-module
  $\Gamma^{4,2}\C^3$, corresponding to the Young diagram
  
  $$\yng(4,2)$$
  
  A presentation of $\sym^2T$ can be described as follows. First we recall that $T$ is presented by the following exact sequence
  
  $$0\rightarrow\Oc\rightarrow\C^3\otimes \Oc(1)\rightarrow T\rightarrow 0$$
 which  dualizes to
  $$0\rightarrow T^\vee\rightarrow\C^3\otimes\Oc(-1)\rightarrow\Oc\rightarrow 0$$
  and since $T^\vee = T(-3)$ we get
  $$0\rightarrow T\rightarrow\C^3\otimes\Oc(2)\rightarrow\Oc(3)\rightarrow 0$$
       
  The first and third sequence fit together into the presentation of $T$
  
  $$ \begin{array}{ccccc}\C^3\otimes\Oc(1)&&\rig{f}&&\C^3\otimes\Oc(2)\\& \searrow&&\nearrow\\ &&T\end{array} $$
  where the horizontal skew-symmetric map $f(v)=v\wedge x$ is given by the matrix
  $\begin{pmatrix}0&x_2&-x_1\\ -x_2&0&x_0\\ x_1&-x_0&0\end{pmatrix}$
  
  The second symmetric power $\sym^2T=E$ appears in the exact sequence
   $$0\rightarrow\C^3\otimes\Oc(1)\rightarrow\sym^2\C^3\otimes \Oc(2)\rightarrow \sym^2T\rightarrow 0$$
  
  and repeating the above argument, thanks to the identification $\sym^2T^\vee = \sym^2T(-6)$ we get the presentation
  
   $$ \begin{array}{ccccc}\sym^2\C^3\otimes\Oc(2)&&\rig{\sym^2f}&&\sym^2\C^3\otimes\Oc(4)\\& \searrow&&\nearrow\\ &&\sym^2T\end{array} $$
  where the horizontal symmetric  map $\sym^2f(v^2)=(v\wedge x)^2$ is given by the matrix
  \begin{equation}\label{eq:B}B=\begin{pmatrix}
      0&0&0&{z}_{2}^{2}&{-2\,{z}_{1}{z}_{2}}&{z}_{1}^{2}\\
      0&{-2\,{z}_{2}^{2}}&2\,{z}_{1}{z}_{2}&0&2\,{z}_{0}{z}_{2}&{-2\,{z}_{0}{z}_{1}}\\
      0&2\,{z}_{1}{z}_{2}&{-2\,{z}_{1}^{2}}&{-2\,{z}_{0}{z}_{2}}&2\,{z}_{0}{z}_{1}&0\\
      {z}_{2}^{2}&0&{-2\,{z}_{0}{z}_{2}}&0&0&{z}_{0}^{2}\\
      {-2\,{z}_{1}{z}_{2}}&2\,{z}_{0}{z}_{2}&2\,{z}_{0}{z}_{1}&0&{-2\,{z}_{0}^{2}}&0\\
      {z}_{1}^{2}&{-2\,{z}_{0}{z}_{1}}&0&{z}_{0}^{2}&0&0\end{pmatrix}\end{equation}

  Note that for $L=\Oc(6)$ we have $E=E^{\vee}\otimes L$ hence, as in \cite{LandOtt13}
  we have a contraction map 
  
   \begin{equation}\label{eq:youngcontraction}\begin{array}{ccc}End(\sym^2\C^3)&\rig{P_f}&End(\sym^2\C^3)\\
  \downarrow&&\uparrow\\
  H^0(E)&\rig{A_f} &H^0(E^\vee\otimes L)^\vee\end{array}\end{equation}
  
  where the horizontal map $P_f$ for $f=v^6\in v_6(\Pj^2)$ is defined as
  $$P_{v^6}(M^2)(w^2)=(M(v)\wedge v\wedge w)^2v^2\quad\forall M\in End(\C^3), w^2\in\sym^2\C^3$$
 where $M^2\in End(\sym^2\C^3)$ is defined by $M^2(v^2)=\left(M(v)\right)^2$ and then extended by linearity
 to any $N\in End(\sym^2\C^3)$ and to any plane sextic $f$.
  Comparing with \cite[\S 2]{Ott09} we see that the invariant $\det A_f$ of degree $27$ of plane sextics  has a construction  similar to the Aronhold invariant
  of degree $4$  of plane cubics.
  
  The coordinate description of $P_f$ is the following. Differentiate  the $6\times 6$ catalecticant $C(f)$ (given by differentiating on 
  rows and columns by monomials of degree $2$) by $B$ in (\ref{eq:B}) .
  The output is a $36\times 36$ matrix obtained by replacing any entry of $B$ by a $6\times 6$ catalecticant block.
  
 The M2 \cite{Macaulay2} commands to get $P_f$ are the following, after the matrix $B$ has been defined as in (\ref{eq:B})
  \begin{verbatim}
  R=QQ[z_0..z_2]
  f=z_0^3+z_1^3+z_2^3+5*z_0*z_1*z_2---any cubic polynomial
  P=diff(transpose basis(2,R),diff(basis(2,R),diff(B,f)))
  \end{verbatim}

  The map $A_f$ is symmetric and could be obtained by a convenient submatrix of $P_f$.
  An alternative way to get a coordinate description is to employ the M2 package ``PieriMaps'' by Steven Sam  with the M2 commands \cite{Macaulay2}
  \begin{verbatim}
loadPackage "PieriMaps"
pieri({6,4,2},{1,1,2,2,3,3},QQ[x_0..x_2]);
  \end{verbatim}
  
  although the symmetry is not transparent from the coordinates chosen by the system.
  
  The contraction $A_f$  in (\ref{eq:youngcontraction}) can be pictorially described (compare with \cite{Ott09}) by

$$\yng(4,2)\otimes\quad\young(******)\rightarrow\young(\ \ \ \ **,\ \ **,**)\quad\simeq\quad \yng(4,2)$$

We have the decomposition
$$End(\sym^2\C^3)=\Gamma^{4,2}\C^3\oplus \Gamma^{2,1}V\oplus\wedge^3\C^3$$
corresponding to
$$\yng(2,2)\otimes\yng(2)=\yng(4,2)\oplus\yng(3,2,1)\oplus\yng(2,2,2)$$
Note that the last two summands correspond to $End(\C^3)$.

If $f\in v_3(\Pj^2)$ then ${\mathrm rk}(A_f)={\mathrm rk}(E)=3$, so that if $f\in S^8(v_3(\Pj^2)$ then ${\mathrm rk}(A_f)\le 3\cdot 8=24$ and in 
particular $\det A_f=0$.  This was observed already in \cite[Theorem 4.2.9]{LandOtt13}.

\begin{proposition}\label{prop:eqW}\begin{enumerate}
\item{} The variety $W$ of sextics in $S^9(v_3(\Pj^2))$ such that their Waring decomposition comes from a complete intersection
of two cubics is cut by the $9$-minors of $\Cat^3$. It has codimension $3$ and degree $165$ and it coincides with the singular locus of $S^9(v_3(\Pj^2))$.
\item{} The secant variety $S^8(v_3(\Pj^2))$ coincides with the reduced structure on $W\cap V(\det A_f)$, it has
codimension $4$ and degree $1485=165\cdot 9$. 
\end{enumerate}
\end{proposition}

\begin{proof} First recall that the variety of symmetric $10\times 10$ matrices of rank $\le 8$ is irreducible
of codimension $3$ and degree $165$ by Segre formula \cite{Segre00}.  Since the variety $\Cat^3_9$ given by the $9$-minors of $\Cat^3$ 
has again codimension $3$, being a linear section of the previous one, it has the same degree $165$. This variety $\Cat^3_9$  corresponds to the
row with $s=8$ of Table 3.1 in \cite{IK}, indeed with the notations of \cite{IK}, it is the union of irreducible strata $Gor(T)$ for several Hilbert functions $T$,
and one may check using Conca-Valla formula ( \cite{ConcaValla99} or \cite[Theorem 4.26]{IK}) that only one strata has the maximum dimension.
It follows that $\Cat^3_9$  is irreducible.
We know that $W$ is irreducible and contained in the variety given by the $9$-minors of $\Cat^3$ ,
since the two varieties have the same dimension, the equality follows, proving i). The claim about the singular locus follows because by direct computation
on a random linear subspace the singular locus of $S^9(v_3(\Pj^2))$, which obviously contains $W$, has codimension $3$ and degree $165$.

In order to prove ii), recall that if $f\in S^8(v_3(\Pj^2))$
then $\det(A_f)=0$. Pick $9$ general points $p_i$ for $i=0,\ldots, 8$ and consider the linear span $f=\sum_{i=0}^8k_ip_i^6$
with coefficients $k_i$, in other words we restrict $f$ to a general $9$-secant space. The expression $\det A_f$ is symmetric in $k_i$,
has total degree $27$ and vanishes when $\prod_{i=0}^8 k_i$ vanishes. An explicit computation with the coordinate expression 
found by the package ``PieriMaps'' as described above shows that $\det(A_f)$ coincides, up to scalar multiples, with
$\left(\prod_{i=0}^8 k_i\right)^3$. We get that the condition that $\det(A_f)$ vanishes on a general $9$-secant space is equivalent to $f$ being of rank $8$.
This expression is unchanged if $f\in W$, namely if $f$ belongs to a $9$-secant space corresponding to a complete intersection of two general cubics.
This computation shows that $V(\det A_f)$ meets $W$ with multiplicity $3$ at a general point of the intersection, so that the reduced structure on 
$W\cap V(\det A_f)$ has
 degree $1485=165\cdot 9$.  Note the rank of $A_f$ drops by $3$ on $W$The degree of $S^8(v_3(\Pj^2))$ can be computed with Numerical Algebraic Geometry.
 The M2 package {\it NumericalImplicitization} by Cho and Kileel shows indeed that its degree is $1485$, which can be certified by the trace test. 
  Since $S^8(v_3(\Pj^2))\subseteq \left(W\cap V(\det A_f)\right)_{red}$, this proves ii)
 \end{proof}

\begin{remark}
ii) of Prop. \ref{prop:eqW} strengthens \cite[Theorem 4.2.9]{LandOtt13} where the $25$-minors of $A_f$ were considered.
\end{remark}

We wish now to find the ramification locus of the $2:1$ map $s^9:AS^9(v_6(\Pj^2))\to S^9(v_6(\Pj^2))$.
Let $p_1,\ldots, p_{10}$ be points in $\Pj^2$.
Let $C(p_1,\ldots, p_{10})$ be a multihomogeneous polynomial of degree $3$ in the coordinates of the points $p_1, . . . , p_{10}$ , skew-symmetric with 
respect to them and that
vanishes if and only if $p_1, . . . , p_{10}$ are on a cubic. A coordinate expression is given by the $10\times 10$ determinant such that
its $i$-th row is the evaluation on $P_i$ of the monomial basis of $\sym^3\C^3$.

We construct now the Terracini matrix of the points $p_1,\ldots, p_9$ which defines the tangent space at $S^9(v_6(\Pj^2)$ at any point in 
$\langle p_1^6,\ldots p_9^6\rangle$.
Consider the $3\times 28$ Jacobian matrix $J$ of the monomial basis  of $\sym^6\C^3$.
Let $J(p_i)$ be the evaluation of $J$ at $p_i$ and let $T(p_1,\ldots, p_9)$ be the $27\times 28$ Terracini matrix obtained by stacking $J(p_1),\ldots, J(p_9)$.
Let $R(p_1,\ldots, p_9;p_{10})=\det\left(\begin{array}{ccc}&T(p_1,\ldots, p_9)\\ \hline \\(p_{10})_0^6&\ldots&(p_{10})_2^6\end{array}\right)$
be the determinant of the $28\times 28$ matrix obtained by stacking $T(p_1,\ldots, p_9)$ and the monomial basis  of $\sym^6\C^3$ evaluated at $p_{10}$.

$R(p_1,\ldots, p_9;p_{10})$ is a multihomogeneous polynomial of degree $15$ in the coordinates of the points 
$p_1, . . . , p_{9}$ , of degree $6$ in the coordinates of $p_{10}$, skew-symmetric with respect to $p_1\ldots, p_9$ and that
vanishes if and only if $p_{10}$ lies on a sextic singular at $p_1, . . . , p_{9}$.

Since $C(p_1,\ldots, p_{10})^2$ is a (non reduced) sextic singular at $p_1,\ldots, p_{10}$ we have the factorization

$$ R(p_1,\ldots, p_9;p_{10}) = C(p_1,\ldots, p_{10})^2N(p_1,\ldots, p_9)$$
where $N(p_1,\ldots, p_9)$  is a multihomogeneous polynomial of degree $9$ in the coor-
dinates of the points $p_1, . . . , p_{9}$ , symmetric with respect to $p_1\ldots, p_9$ and that
vanishes if and only if there is a reduced sextic  singular at $p_1, . . . , p_{9}$. Note that the result in \cite{ArbarelloCornalba81}
that the unique sextic singular at $p_1,\ldots, p_9$ is the non-reduced curve $C(p_1,\ldots, p_9,p)^2$ in $p$ is equivalent to the fact that $T$ has maximal rank
for a general choice of $p_1,\ldots, p_9$, which can be checked by a random choice. This proves that $N(p_1,\ldots, p_9)$ is a nonzero polynomial,
which defines the codimension $1$ condition that the space of sextics singular at $p_1,\ldots, p_9$ ha dimension $\ge 2$.
It can be checked that at general  $p_1,\ldots, p_9$ on the hypersurface  $N(p_1,\ldots, p_9)=0$ then $T$ has corank $1$, while if 
$p_1,\ldots, p_9$ are distinct points defined by a general complete intersection
then both  $R(p_1,\ldots, p_9;p_{10})$ and $C(p_1,\ldots, p_{10})$ vanish $\forall P_{10}$ and $N(p_1,\ldots, p_9)\neq 0$.

We note also that the same argument given in the proof of \cite[Theor. 45]{OttSernesi10} shows that $N(p_1,\ldots, p_9)$ is irreducible, since there are no 
nonzero symmetric cubic invariants of $9$ points in $\Pj^2$.

\begin{remark} For general $p_1,\ldots, p_8$, the locus of ninth point $p_9$ such that there is a reduced sextic singular at $p_1,\ldots, p_8, p_9$
has two irreducible components, namely the nonic $N(p_1,\ldots, p_9)=0$ and the point $p_9$ in the base locus of the linear system
of cubics through $p_1,\ldots, p_8$. The general sextic arising from the first (resp. second) component is irreducible (resp. reducible)
and lies in $\Rm$ (resp. in $W$), see Definition \ref{def:RW}. More precisely, $\Rm\subset S^9(v_6(\Pj^2))$ is the closure of $\cup \langle p_1^6,\ldots, p_9^6\rangle$, 
where the union is taken for distinct $p_1,\ldots, p_9$ 
satisfying $N(p_1,\ldots, p_9)=0$. As noted in Remark \ref{RW}, $\Rm$ and $W$ are two irreducible components of the ramification locus of the $2:1$ map $AS^9(v_6(\Pj^2))\to S^9(v_6(\Pj^2))$. They are distinct since $H_{27}$ vanishes on $\Rm$ and not on $W$.
When $N(p_1,\ldots, p_9)=0$ then $Z=\{p_1,\ldots, p_9\}$ is self-linked on the cubic $C(p_1,\ldots, p_9,p)$ with respect to any irreducible sextic
singular at $p_1,\ldots, p_9$, which meets the cubic in $2p_1+\ldots +2p_9$.
\end{remark}

\begin{theorem}
Let $f=\sum_{i=1}^9p_i^6$ be a rank $9$ sextic, in affine notation.
Then $$\det A_f=\lambda N(p_1,\ldots p_9)^2$$ for a nonzero scalar $\lambda\in\C^*$.
\end{theorem}

\begin{proof} Recall we denoted $E=\sym^2T$ and denote $Z=\{p_1,\ldots, p_9\}$. Assume that $\det A_f$ vanishes, then $\ker A_f$ is a nonzero sibspace of
$ H^0(E)$.
Consider the restriction $H^0(E)\rig{j} H^0(E_{|Z})$, note both sides are $27$-dimensional. Assume $j$ is injective. Then by  \cite[Lemma 5.4.1]{LandOtt13} we have that
$H^0(I_Z\otimes E)=\ker A_f$ is nonzero, which is a contradiction since a section vanishing on $Z$ is in the kernel of $j$.
It follows that $j$ is not injective and there is a section $s\in H^0(E)$ vanishing at $Z=\{p_1,\ldots, p_9\}$. Again by \cite[Lemma 5.4.1]{LandOtt13} 
the section $s$ belongs to $\ker A_f$. Then the contraction $$\sym^2 H^0(I_Z\otimes E)\to H^0(I_Z\otimes E)\otimes H^0(I_Z\otimes E^\vee\otimes\Oc(6))\to H^0(I_{Z^2}\otimes\Oc(6))$$ as in \cite[Theorem 5.4.3]{LandOtt13}
takes $s^2$ to a sextic singular at $Z$. There are two cases, $Z$ is a complete intersection or $N(p_1,\ldots, p_9)=0$. The first case 
has codimension $3$ by \ref{prop:eqW} and can be excluded in proving a polynomial equality, like in our statement.
It follows that $\det A_f=0$ implies generically $N(p_1,\ldots, p_9)=0$, hence a power of $ N(p_1,\ldots p_9)$ is divided by $\det A_f$.
The multihomogeneous polynomial $\det A_f$ is symmetric in $p_1,\ldots, p_9$ and has total degree $27\cdot 6$, hence must have degree $18$ in each $p_i$.
The result follows since $N(p_1,\ldots, p_9)$ is irreducible.
\end{proof}

\begin{corollary}
The variety $\Rm$ is the complete intersection of the two hypersurfaces $\det A_f$ and $S^9(v_6(\Pj^2))$ and has degree $10\cdot 27= 270 $.
\end{corollary}

\vskip 1cm

{\bf Tables about the loci we have studied}

 In the next table we list some informations about the loci we have studied in $\Pj^{27}=\Pj(\sym^6\C^3)$. The main theme is that $(\Cat)_3$ and the invariant $H_{27}=\det A_f$ are enough to describe
all $k$-secant varieties to $v_6(\Pj^2)$.

$$\begin{array}{c|c|c|c|c}
&dim&deg&\textrm{equations}\\
\hline\\
S^9&26&10&\det\Cat^3&\textrm{9-secant}\\
\Rm&25&270&\det\Cat^3, H_{27}&\textrm{comp. of ramif. locus of\  }AS^9\to S^9\\
W&24&165&\left(\Cat^3\right)_9&\mathrm{Sing\ } S^9, \textrm{span of compl. inters. 9-ples}\\
S^8&23&1485=(165\cdot 27)/3&\left(\Cat^3\right)_9, H_{27}&\textrm{8-secant}\\
W'&21&2640&\left(\Cat^3\right)_8&\textrm{sextics with three apolar cubics}\\
S^7&20&11880=(2640\cdot 27)/6&\left(\Cat^3\right)_8, H_{27}&\textrm {7-secant}\end{array}$$

The fact that $\mathrm{Sing\ } S^9=W$ has been checked computationally,
since the singular locus of $S^9$ has codimension $3$ and degree $165$, and contains the variety $W$ given by the $9$-minors of $\Cat^3$,
hence the equality follows. The following graph links the containments among the several loci.

$$\begin{array}{ccc}S^9\\
\vert&\diagdown\\
W&&\Rm\\
\vert&\diagdown&\vert\\
W'&&S^8\\
&\diagdown&\vert\\
&&S^7\end{array}$$

The right colums is obtained by the left column after cut with the invariant hypersurface $H_{27}$. The three diagonal links correspond to cut with the invariant $H_{27}$. The matrix $A_f$ drops rank by $1$ in the first link,
by $3$ in the second one, by $6$ in the third one. Hence a check on the degrees shows that the intersection with $H_{27}$ does not contain other components
in all the three cases.
\vskip 0.8cm
The following table lists the properties related to Waring decomposition of ge\-neral member of the loci ($\mathrm{ brk}$ is the border rank)
$$\begin{array}{c|c|c|c|l}
&\mathrm{ rk\ }\Cat^3&\mathrm {brk}&\mathrm {rk}&\#\{\textrm{ minimal Waring decompositions}\}\\
\hline\\
S^9&9&9&9&2\\
\Rm&9&9&9&1\\
W&8&9&9&1\textrm{\ (Prop.\ref{2cub3}) } \\
S^8&8&8&8&1\\
W'&7&9&9&\infty\\
S^7&7&7&7&1\end{array}$$
\vskip 0.8cm
{\bf Tables about the dual loci}

$$\begin{array}{c|c|c|c}
&dim&deg&\\
\hline\\
(S^7)^\vee&20&34435125\textrm{\cite{CapHar98}}&\textrm{irred. sextics with 7 singular pts}\\
(W')^\vee&26&83200\textrm{\cite{BHORS12}}&\textrm{sums of three squares}\\
(S^8)^\vee&19&58444767\textrm{\cite{CapHar98}}&\textrm{irred. sextics with 8 singular pts}\\
W^\vee&18&\frac{1}{2}{{18}\choose 9}=24310&\textrm{reducible in 2 cubics, or sums of two squares}\\
\Rm^\vee&18&57435240\textrm{\cite{Getz97,CapHar98}}&\textrm{irred. sextics with 9 singular pts}\\
(S^9)^\vee&9&2^9=512&\textrm{square of a cubic}\end{array}
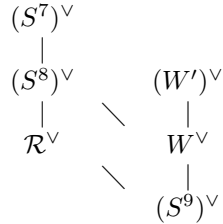
$$

$$\begin{array}{ccc}(S^7)^\vee\\
\vert&\\
(S^8)^\vee&&(W')^\vee\\
\vert&\diagdown&\vert\\
\Rm^\vee&&W^\vee\\
&\diagdown&\vert\\
&&(S^9)^\vee\end{array}$$

\bibliographystyle{amsplain}
\bibliography{biblioLuca}

\providecommand{\bysame}{\leavevmode\hbox to3em{\hrulefill}\thinspace}
\providecommand{\MR}{\relax\ifhmode\unskip\space\fi MR }
\providecommand{\MRhref}[2]{%
  \href{http://www.ams.org/mathscinet-getitem?mr=#1}{#2}
}
\providecommand{\href}[2]{#2}
\begin{thebibliography}{10}

\bibitem{AlexHir95}
J.~Alexander and A.~Hirschowitz, \emph{Polynomial interpolation in several
  variables}, J. Algebraic Geom. \textbf{4} (1995), 201--222.

\bibitem{AngeC20}
E.~Angelini and L.~Chiantini, \emph{On the identifiability of ternary forms},
  Lin. Alg. Applic. \textbf{599} (2020), 36--65.

\bibitem{AngeCMazzon19}
E.~Angelini, L.~Chiantini, and A.~Mazzon, \emph{Identifiability for a class of
  symmetric tensors}, Mediterr. J. Math. \textbf{16} (2019), 97.

\bibitem{ArbarelloCornalba81}
E.~Arbarello and M.~Cornalba, \emph{Footnotes to a paper of {B}. {S}egre},
  Math. Ann. \textbf{256} (1981), 341--362.

\bibitem{Ball19}
E.~Ballico, \emph{An effective criterion for the additive decompositions of
  forms}, Rend. Ist. Matem. Trieste \textbf{51} (2019), 1--12.

\bibitem{BallC}
E.~Ballico and L.~Chiantini, \emph{On the terracini locus of projective
  varieties}, available online arXiv:2011.13189, 2020.

\bibitem{BHORS12}
G.~Blekherman, J.~Hauenstein, J.~C. Ottem, K.~Ranestad, and B.~Sturmfels,
  \emph{Algebraic boundaries of {H}ilbert's {SOS} cones}, Compos. Math.
  \textbf{148} (2012), no.~6, 1717--1735. \MR{2999301}

\bibitem{BucaBucTeit13}
W.~Buczy\'{n}ska, J.~Buczy\'{n}ski, and Z.~Teitler, \emph{Waring decompositions
  of monomials}, J. of Algebra \textbf{378} (2013), 45--57.

\bibitem{BucGinenskyLand13}
J.~Buczy\'{n}ski, A.~Ginensky, and J.M. Landsberg, \emph{Determinantal
  equations for secant varieties and the {E}isenbud-{K}oh-{S}tillman
  conjecture}, J. London Math. Soc. \textbf{88} (2013), 1--24.

\bibitem{CapHar98}
L.~Caporaso and J.~Harris, \emph{Counting plane curves of any genus}, Invent.
  Math. \textbf{131} (1998), no.~2, 345--392. \MR{1608583}

\bibitem{C19}
L.~Chiantini, \emph{Hilbert functions and tensor analysis}, Quantum Physics and
  Geometry, Lecture Notes of the Unione Matematica Italiana, vol.~25, Springer,
  Berlin, New York NY, 2019, pp.~125--151.

\bibitem{CCi02a}
L.~Chiantini and C.~Ciliberto, \emph{Weakly defective varieties}, Trans. Amer.
  Math. Soc. \textbf{354} (2002), 151--178.

\bibitem{CCi06}
\bysame, \emph{On the concept of k-secant order of a variety}, J. London Math.
  Soc. \textbf{73} (2006), 436--454.

\bibitem{COttVan17a}
L.~Chiantini, G.~Ottaviani, and N.~Vannieuwenhoven, \emph{On generic
  identifiability of symmetric tensors of subgeneric rank}, Trans. Amer. Math.
  Soc. \textbf{369} (2017), 4021--4042.

\bibitem{ConcaValla99}
A.~Conca and G.~Valla, \emph{Hilbert function of powers of ideals of low
  codimension}, Math. Zeit. \textbf{230} (1999), 753--784.

\bibitem{Davis84}
E.~Davis, \emph{Hilbert functions and complete intersections}, Rend. Seminario
  Mat. Univ. Politecnico Torino \textbf{42} (1984), 333--353.

\bibitem{Getz97}
E.~Getzler, \emph{Intersection theory on {$\overline{ M}_{1,4}$} and elliptic
  {G}romov-{W}itten invariants}, J. Amer. Math. Soc. \textbf{10} (1997), no.~4,
  973--998. \MR{1451505}

\bibitem{Macaulay2}
D.~Grayson and M.~Stillman, \emph{Macaulay 2, a software system for research in
  algebraic geometry}, Available at \/ {\tt http://www.math.uiuc.edu/Macaulay2/
  }.

\bibitem{Harr86}
J.~Harris, \emph{On the {S}everi problem}, Invent. Math. \textbf{84} (1986),
  445--461.

\bibitem{IK}
A.~Iarrobino and V.~Kanev, \emph{Power sums, {G}orenstein algebras, and
  determinantal loci}, Lecture Notes in Mathematics, vol. 1721, Springer,
  Berlin, New York NY, 1999.

\bibitem{LandOtt13}
J.M. Landsberg and G.~Ottaviani, \emph{Equations for secant varieties of
  {V}eronese and other varieties}, Ann. Mat. Pura Appl. \textbf{192} (2013),
  569--606.

\bibitem{Migliore}
J.~Migliore, \emph{Introduction to liaison theory and deficiency modules},
  Progress in Mathematics, vol. 165, Birk{\"a}user, Basel, Boston MA, 1998.

\bibitem{MourOneto20}
B.~Mourrain and A.~Oneto, \emph{On minimal decompositions of low rank symmetric
  tensors}, Lin. Alg. Applic. \textbf{607} (2020), 347--377.

\bibitem{OedOtt13}
L.~Oeding and G.~Ottaviani, \emph{Eigenvectors of tensors and algorithms for
  {W}aring decomposition}, J. Symbolic Comput. \textbf{54} (2013), 9--35.

\bibitem{Ott09}
G.~Ottaviani, \emph{An invariant regarding {W}aring's problem for cubic
  polynomials}, Nagoya Math. J. \textbf{193} (2009), 95--110.

\bibitem{OttSernesi10}
G.~Ottaviani and E.~Sernesi, \emph{On the hypersurface of {L}\"{u}roth
  quartics}, Michigan Math. J. \textbf{59} (2010), no.~2, 365--394.

\bibitem{Segre00}
C.~Segre, \emph{Gli ordini delle variet{\`a} che annullano i determinanti dei
  diversi gradi estratti da una data matrice}, Atti Accad. Lincei Classe Sci.
  \textbf{9} (1900), 253--260.

\end{thebibliography}

\end{document}